\documentclass[a4paper,10pt]{article}
  \usepackage[light]{antpolt} 
  \usepackage{apacite} 
  \usepackage[T1]{fontenc}
  \usepackage{sectsty} 
  \usepackage{hyperref} 
  \usepackage{graphicx} 
  \usepackage{comment} 
  \usepackage{twoopt} 
  \usepackage{setspace} 
  \usepackage{enumitem} 
  \usepackage{amsbsy} 
  \usepackage[sectionbib]{natbib} 
  \usepackage{todonotes} 
  \usepackage{wrapfig} 
  \usepackage{tocloft} 
  \usepackage{float} 

  \usepackage{amsmath, amssymb, amsfonts} 
  \usepackage{amsthm} 
  \usepackage{mathrsfs}	
  \usepackage{stmaryrd}	

  \usepackage{wallpaper} 
  \usepackage{fancyhdr}	
  \usepackage{changepage} 

  \usepackage{todonotes} 

\setlist{nolistsep}

\newtheoremstyle{scthmstyle} 
	{15pt} 
	{15pt} 
	{\itshape} 
	{} 
	{\bfseries\scshape} 
	{.} 
	{.5em} 
	{} 

\newtheoremstyle{scdefstyle} 
	{15pt} 
	{15pt} 
	{\normalfont} 
	{} 
	{\bfseries\scshape} 
	{.} 
	{.5em} 
	{} 

\newtheoremstyle{scremstyle} 
	{15pt} 
	{15pt} 
	{\normalfont} 
	{} 
	{\itshape} 
	{.} 
	{.5em} 
	{} 

\newtheoremstyle{scclaistyle} 
	{15pt} 
	{15pt} 
	{\normalfont} 
	{0.5cm} 
	{\itshape} 
	{.} 
	{.5em} 
	{} 

\theoremstyle{scthmstyle}
\newtheorem{theorem}{Theorem}[section]
\newtheorem{proposition}[theorem]{Proposition}
\newtheorem{lemma}[theorem]{Lemma}
\newtheorem{corollary}[theorem]{Corollary}
\theoremstyle{scdefstyle}
\newtheorem{definition}[theorem]{Definition}

\newtheorem{example}[theorem]{Example}
\newtheorem{question}[theorem]{Question}
\theoremstyle{scremstyle}
\newtheorem{remark}[theorem]{Remark}
\theoremstyle{scclaistyle}
\newtheorem{claim}{Claim}[theorem]

  \newcommand{\eq}[1]{\begin{align*} #1 \end{align*}}

  \renewcommand{\abstract}[1]{\begin{quote}{\footnotesize\textsc{Abstract.} #1}\\\end{quote}}

  \newcommandtwoopt{\theo}[3][][]{
    \begin{theorem}[#1]\label{#2}
      #3
    \end{theorem}}
  \newcommandtwoopt{\prop}[3][][]{
    \begin{proposition}[#1]\label{#2}
      #3
    \end{proposition}}
  \newcommandtwoopt{\lemm}[3][][]{
    \begin{lemma}[#1]\label{#2}
      #3
    \end{lemma}}
  \newcommandtwoopt{\coro}[3][][]{
    \begin{corollary}[#1]\label{#2}
      #3
    \end{corollary}}
  \newcommandtwoopt{\defi}[3][][]{
    \begin{definition}[#1]\label{#2}
      #3$\hfill\dashv$
    \end{definition}}
  \newcommandtwoopt{\exam}[3][][]{
    \begin{example}[#1]\label{#2}
      #3
    \end{example}}
  \newcommandtwoopt{\ques}[3][][]{
    \begin{question}[#1]\label{#2}
      #3
    \end{question}}
  \newcommandtwoopt{\rema}[3][][]{
    \begin{remark}[#1]\label{#2}
      #3
    \end{remark}}

  \newcommandtwoopt{\qtheo}[3][][]{
    \begin{theorem}[#1]\label{#2}
      #3$\hfill\dashv$
    \end{theorem}}
  \newcommandtwoopt{\qprop}[3][][]{
    \begin{proposition}[#1]\label{#2}
      #3$\hfill\dashv$
    \end{proposition}}
  \newcommandtwoopt{\qlemm}[3][][]{
    \begin{lemma}[#1]\label{#2}
      #3$\hfill\dashv$
    \end{lemma}}
  \newcommandtwoopt{\qcoro}[3][][]{
    \begin{corollary}[#1]\label{#2}
      #3$\hfill\dashv$
    \end{corollary}}

  \newcommandtwoopt{\defin}[3][][]{
    \begin{definition}[#1]\label{#2}
      #3
    \end{definition}}

  \renewcommand{\proof}[1]{\textsc{Proof.} #1$\qed$\\}

  \renewcommand{\qed}{\hfill\blacksquare}

  \DeclareMathOperator{\M}{\mathcal M}
  \DeclareMathOperator{\N}{\mathcal N}

  \DeclareMathOperator{\C}{\mathcal C}

  \DeclareMathOperator{\ran}{ran}

  \DeclareMathOperator{\con}{Con}

  \DeclareMathOperator{\crit}{crit}

  \DeclareMathOperator{\contr}{\lightning}

  \DeclareMathOperator{\forces}{\Vdash}

  \DeclareMathOperator{\restr}{\upharpoonright}

  \DeclareMathOperator{\col}{Col}

  \renewcommand{\P}{\mathcal P}

  \renewcommand{\subset}{\subseteq}

	\newcommand{\p}{\mathscr P}

  \newcommand{\bra}[1]{\langle #1\rangle}

  \newcommand{\abs}[1]{\left|#1\right|}

  \newcommand{\pnormal}{\mathrel{\ooalign{$\lneq$\cr\raise.22ex\hbox{$\lhd$}\cr}}}
  \newcommand{\pideal}{\mathrel{\ooalign{$\lneq$\cr\raise.22ex\hbox{$\lhd$}\cr}}}

  \newcommand{\zf}{\textsf{ZF}}
  \newcommand{\zfc}{\textsf{ZFC}}

  \newcommand{\on}{\textsf{On}}

  \newcommand{\gbc}{\textsf{GBC}}

  \newcommand{\dc}{\textsf{DC}}

  \newcommand{\godel}[1]{\ulcorner #1 \urcorner}

  \newcommand{\po}{\ar@{}[dr]|{\text{\pigpenfont R}}}
  \newcommand{\pb}{\ar@{}[dr]|{\text{\pigpenfont J}}}

\newcommand\Vopenka{Vop\v{e}nka}
\begin{document}

\title{
  \vspace{-3cm}
  The Virtual Large Cardinal Hierarchy
}
\author{
    Stamatis Dimopoulos$^1$,
    Victoria Gitman$^2$ and
    Dan Saattrup Nielsen$^3$
}

\date{
  \begin{small}
    $^1$ Research and Innovation Foundation, Nicosia, Cyprus\\
    $^2$ CUNY Graduate Center, New York, USA\\
    $^3$ The Alexandra Institute, Copenhagen, Denmark\\[2ex]
  \end{small}
  \today\thanks{We would like to thank the anonymous referee for their careful reading of the paper, excellent suggestions for improvements, and new arguments.}
}
\maketitle

\abstract{
We continue the study of the virtual large cardinal hierarchy by analysing virtual
versions of superstrong, Woodin, and Berkeley cardinals. Gitman and Schindler showed
that virtualizations of strong and supercompact cardinals yield the same large cardinal
notion. We provide various equivalent characterizations of virtually Woodin cardinals,
including showing that $\on$ is virtually Woodin if and only if for every class $A$,
there is a proper class of virtually $A$-extendible cardinals. We introduce the virtual
\Vopenka\ principle for finite languages and show that it is not equivalent to the
virtual \Vopenka\ principle (although the two principles are equiconsistent), but is
equivalent to the assertion that $\on$ is virtually pre-Woodin, a weakening of
virtually Woodin, which is equivalent to having for every class $A$, a weakly virtually
$A$-extendible cardinal. We show that if there are no virtually Berkeley cardinals,
then $\on$ is virtually Woodin if and only if $\on$ is virtually pre-Woodin (if and
only if the virtual \Vopenka\ principle for finite languages holds). In particular, if
the virtual \Vopenka\ principle holds and $\on$ is not Mahlo, then $\on$ is not
virtually Woodin, and hence there is a virtually Berkeley cardinal.
}

\begin{spacing}{0.0}
  \tableofcontents
\end{spacing}

\section{Introduction}

The study of generic large cardinals, being cardinals that are critical points of
elementary embeddings existing in generic extensions, goes back to the 1970's. At that
time, the primary interest was the existence of \textit{precipitous} and
\textit{saturated} ideals on small cardinals like $\omega_1$ and $\omega_2$. Research
in this area later moved to the study of more general generic embeddings, both defined
on $V$, but also on rank-initial segments of $V$ --- these were investigated by e.g.
\cite{DonderLevinski} and \cite{FerberGitik}.

\qquad The move to \textit{virtual} large cardinals happened when \cite{Schindler}
introduced the \textit{remarkable cardinals}, which it turned out later were precisely
a virtualization of supercompactness. Various other virtual large cardinals were first
investigated in \cite{GitmanSchindler}. The key difference between virtual large
cardinals and generic versions of large cardinals studied earlier is that in the
virtual case we require the embedding to be between sets with the target model being a
subset of the ground model. These assumptions imply that virtual large cardinals are
actual large cardinals: they are at least ineffable, but small enough to exist in $L$.
These large cardinals are special because they allow us to work with embeddings as in
the higher reaches of the large cardinal hierarchy while being consistent with $V=L$,
which enables equiconsistencies at these ``lower levels''.

\qquad To take a few examples, \cite{Schindler} has shown that the existence of a
remarkable cardinal is equiconsistent with the statement that the theory of $L(\mathbb
R)$ cannot be changed by proper forcing, which was improved to semi-proper forcing in
\cite{Schindler2}. \cite{Wilson} has shown that the existence of a virtually \Vopenka\
cardinal is equiconsistent with the hypothesis
\eq{
  \zf + \godel{\text{${\bf\Sigma}^1_2$ is the class of all $\omega_1$-Suslin sets}} +
  \Theta = \omega_2,
}

and \cite{SchindlerWilson} has shown that the existence of a virtually Shelah for
supercompactness cardinal is equiconsistent with the hypothesis
\eq{
  \zfc + \godel{\text{every universally Baire set of reals has the perfect set property}}.
}

Kunen's Inconsistency fails for virtual large cardinals in the sense that a forcing
extension can have elementary embeddings $j:V_\alpha^V\to V_\alpha^V$ with $\alpha$
much larger than the supremum of the critical sequence. In the theory of large
cardinals, Kunen's Inconsistency is for instance used to prove that requiring that
$j(\kappa)>\lambda$ in the definition of $\kappa$ being $\lambda$-strong is
superfluous. It turns out that the use of Kunen's Inconsistency in that argument is
actually essential because versions of virtual strongness with and without that
condition are not equivalent; see e.g. Corollary~\ref{coro.kunen-prestrong-strong}. The
same holds for virtual versions of other large cardinals where this condition is used
in the embeddings characterization. Each of these virtual large cardinals therefore has
two non-equivalent versions, with and without the condition.

\qquad In this paper, we continue the study of virtual versions of various large
cardinals.

\qquad In Section~\ref{sec:supercompact} we establish some new relationships between
virtual large cardinals that were previously studied. We prove the Gitman-Schindler
result, alluded to in \cite{GitmanSchindler}, that virtualizations of strong and
supercompact cardinals are equivalent to remarkability. We show how the existence of
virtually supercompact cardinals without the $j(\kappa)>\lambda$ condition is related
to the existence of virtually rank-into-rank cardinals.

\qquad In Section \ref{sect.woodvop} we study virtual versions of Woodin cardinals and
introduce the virtual \Vopenka\ principle for finite languages. We provide various
equivalent characterizations of virtually Woodin cardinals. It follows, from the
equivalences, that $\on$ is virtually Woodin if and only if for every class $A$ there
is a virtually $A$-extendible cardinal (as defined in \cite{GitmanHamkins}),
equivalently, there is a stationary class of virtually $A$-extendible cardinals. Recall
from \cite{GitmanHamkins} that the virtual \Vopenka\ Principle holds if and only if for
every class $A$ there is a proper class of weakly virtually $A$-extendible cardinals.
It follows from arguments in \cite{GitmanHamkins} that the virtual \Vopenka\ Principle
for finite languages holds if and only if for every class $A$ there is a weakly
virtually $A$-extendible cardinal. We show that $\on$ is virtually \Vopenka\ for finite
languages if and only if $\on$ is \textit{faintly pre-Woodin}, a weakening of the
notion of virtual Woodinness --- see Definition \ref{defi.woodin}.

\qquad In Section \ref{sect.berkeley} we study a virtual version of \textit{Berkeley}
cardinals, a large cardinal known to be inconsistent with \zfc.  We show that if there
are no virtually Berkeley cardinals, then $\on$ is virtually Woodin if and only if
$\on$ is faintly pre-Woodin if and only if the virtual \Vopenka\ principle for finite
languages holds. In this situation, the virtual \Vopenka\ Principle for finite
languages is equivalent to the virtual \Vopenka\ Principle. However, we will use
virtual Berkeley cardinals to separate the two principles. It follows also that if the
virtual \Vopenka\ Principle holds and $\on$ is not Mahlo (in particular, $\on$ is not
virtually Woodin), then there is a virtually Berkeley cardinal, but as pointed out by
the anonymous referee, it is possible to have that the virtual \Vopenka\ Principle
holds and $\on$ is Mahlo, but $\on$ is not virtually Woodin.

\section{Preliminaries}

We will denote the class of ordinals by $\on$. For sets $X$ and $Y$ we denote by
${^X}Y$ the set of all functions from $X$ to $Y$. For an infinite cardinal $\kappa$, we
let $H_\kappa$ be the set of sets $X$ such that the cardinality of the transitive
closure of $X$ has size less than $\kappa$. The symbol $\contr$ will denote a
contradiction and $\p(X)$ will denote the power set of $X$. We will say that an
elementary embedding $j:\M\to \N$ is \textbf{generic} if it exists in a forcing
extension of $V$. We will use $H_\lambda$ and $V_\lambda$ to denote these sets as
defined in the ground model $V$, while the $H_\lambda$ or $V_\lambda$ of any other
universe will have a superscript indicating which universe it comes from.

\qquad A key folklore lemma which we will frequently need when dealing with generic
elementary embeddings is the following.

\lemm[Countable Embedding Absoluteness][lemm.ctblabs]{
  Let $\M$ and $\N$ be transitive sets and assume that $\M$ is countable. Let
  $\pi\colon\M\to\N$ be an elementary embedding, $\P$ a transitive class with
  $\M,\N\in\P$ and \hbox{$\P\models\zf^-+\dc+\godel{\text{$\M$ is
  countable}}$}.\footnote{The theory $\zf^-$ consists of the axioms of $\zf$ without
  the powerset axiom and with the collections scheme instead of the replacement
  scheme.} Then $\P$ has an elementary embedding $\pi^*\colon\M\to\N$ which agrees with
  $\pi$ on any desired finite set and on the critical point if it exists.
}

The following proposition is an almost immediate corollary of Countable Embedding
Absoluteness.

\prop[][prop.coll]{
Let $\M$ and $\N$ be transitive models and assume that there is a generic elementary
embedding $\pi:\M\to \N$. Then $V^{\col(\omega,\M)}$ has an elementary embedding
$\pi^*:\M\to\N$ which agrees with $\pi$ on any desired finite set and has the same
critical point if it exists.
}

For proofs, see \cite{GitmanSchindler} (Section~3).

\section{Virtually supercompact cardinals}\label{sec:supercompact}
\label{sect.strongsc}
In this section, we establish some relationships between virtual large cardinals
related to virtually supercompact cardinals. We start with definitions of the relevant
virtual large cardinal notions.

\defi{
  Let $\theta$ be a regular uncountable cardinal. Then a cardinal $\kappa<\theta$ is
  \begin{itemize}
    \item \textbf{faintly $\theta$-measurable} if, in a forcing extension, there is a
        transitive set $\N$ and an elementary embedding $\pi\colon H_\theta\to\N$ with
        $\crit\pi=\kappa$,
    \item \textbf{faintly $\theta$-strong} if it is faintly $\theta$-measurable,
        $H_\theta=H_\theta^{\N}$ and $\pi(\kappa)>\theta$,
    \item \textbf{faintly $\theta$-supercompact} if it is faintly $\theta$-measurable,
        ${^{<\theta}\N}\cap V\subset\N$ and $\pi(\kappa)>\theta$.
  \end{itemize}

  We further replace ``faintly'' by \textbf{virtually} when $\N\subset V$, we attach a
  \textbf{``pre''} if we leave out the assumption $\pi(\kappa)>\theta$, and when we do
  not mention $\theta$ we mean that it holds for all regular $\theta>\kappa$. For
  instance, a faintly pre-strong cardinal is a cardinal $\kappa$ such that for all
  regular $\theta>\kappa$, $\kappa$ is faintly $\theta$-measurable with
  $H_\theta\subset\N$.
}

\qquad Observe that whenever we have a virtual large cardinal that has its defining
property for all regular $\theta$, we can assume that the target of the embedding is an
element of the ground model $V$ and not just a subset of $V$. Suppose, for instance,
that $\kappa$ is virtually measurable and fix a regular $\theta>\kappa$ and set
$\lambda:=(2^{<\theta})^+$. Take a generic elementary embedding $\pi:H_\lambda\to
\M_\lambda$ witnessing that $\kappa$ is virtually $\lambda$-measurable. The restriction
$\pi\restr H_\theta\colon H_\theta\to\pi(H_\theta)$ witnesses that $\kappa$ is
virtually $\theta$-measurable and the target model $\M_\theta := \pi(H_\theta)$ is in
$V$ because $\M_\lambda\subseteq V$ by assumption. Thus, the weaker assumption that the
target model $\M_\theta\subseteq V$ only affects level-by-level virtual large
cardinals. Indeed, as we will see in later sections, even further weakening the
assumption $\N\subseteq V$ to $H_\theta=H_\theta^{\N}$ in the definition of virtually
strong (or supercompact) cardinals yields the same notion (again, we do not know
whether this holds level-by-level).

\qquad Small cardinals such as $\omega_1$ can be generically measurable and hence
faintly measurable. However, virtual large cardinals are large cardinals in the usual
sense, as the following shows.

\qquad Recall from \cite{gitman:welch} that a cardinal $\kappa$ is \textbf{1-iterable}
if for every $A\subset\kappa$ there is a transitive $\M\models\zfc^-$ with
$\kappa,A\in\M$ and a weakly amenable $\M$-ultrafilter $\mu$ on $\kappa$ with a
well-founded ultrapower.  Recall that $\mu$ is an \textbf{$\M$-ultrafilter} on $\kappa$
if the structure $\langle M,\in,\mu\rangle$ satisfies that $\mu$ is a normal
ultrafilter on $\kappa$, and such a $\mu$ is \textbf{weakly amenable} if $\mu\cap
X\in\M$ for every $X\in\M$ of $\M$-cardinality $\leq\kappa$. It is not difficult to see
that an $\M$-ultrafilter $\mu$ on $\kappa$ with a well-founded ultrapower is weakly
amenable if and only if the ultrapower embedding $j:\M\to \N$ is
\textbf{$\kappa$-powerset preserving}, meaning that $\M$ and $\N$ have the same subsets
of $\kappa$. 1-iterable cardinals are weakly ineffable limits of ineffable cardinals,
and hence, in particular, weakly compact \cite{gitman:welch}.

\prop[][prop.virtit]{
  For any regular uncountable cardinal $\theta$, every virtually $\theta$-measurable
  cardinal is a 1-iterable limit of 1-iterable cardinals.
}

The proof is essentially the same as the proof of Theorem~4.8 in
\cite{GitmanSchindler}.

\qquad Schindler showed in \cite{schindler:remarkable} that remarkable cardinals can be
viewed as a version of  virtual supercompact cardinals via Magidor's characterization
of supercompactness. Later Gitman and Schindler showed in \cite{GitmanSchindler} that
remarkables are precisely the virtually supercompacts, and indeed surprisingly, they
are also precisely the virtually strongs. So, in particular, virtually strong and
virtually supercompact cardinals are equivalent. We give the proof of the equivalences,
which was omitted in \cite{GitmanSchindler}, here.

\defi{
  Let $\theta$ be a regular uncountable cardinal. Then a cardinal $\kappa<\theta$ is
  \textbf{virtually $\theta$-supercompact ala Magidor} if there are
  $\bar\kappa<\bar\theta<\kappa$ and a generic elementary embedding $\pi\colon
  H_{\bar\theta}\to H_\theta$ such that $\crit\pi=\bar\kappa$ and
  $\pi(\bar\kappa)=\kappa$.
}

\theo[G.-Schindler][theo.rem]{
  For an uncountable cardinal $\kappa$, the following are equivalent.
  \begin{enumerate}
   \item $\kappa$ is faintly strong.
    \item $\kappa$ is virtually strong.
    \item $\kappa$ is virtually supercompact.
    \item $\kappa$ is virtually supercompact ala Magidor.
  \end{enumerate}
}
\proof{
$(iii)\Rightarrow (ii)\Rightarrow(i)$ is simply by definition.

\qquad $(i)\Rightarrow (iv)$: Fix a regular uncountable $\theta>\kappa$ and let
$\delta=(2^{<\theta})^+$. By $(i)$ there exists a generic elementary embedding
$\pi\colon H_\delta\to\M$ with $\crit\pi=\kappa$, $\pi(\kappa)>\delta$, and
\hbox{$H_\delta=H_\delta^{\M}$}. We can restrict the embedding $\pi$ to
$\pi:H_\theta\to H_{\pi(\theta)}^{\M}$. Since $H_\theta, H_{\pi(\theta)}^{\M}\in\M$,
Countable Embedding Absoluteness~\ref{lemm.ctblabs} implies that $\M$ has a generic
elementary embedding $\pi^*:H_\theta\to H_{\pi(\theta)}^{\M}$ with $\crit\pi^*=\kappa$
and $\pi^*(\kappa)=\pi(\kappa)>\theta$. Elementarity of $\pi$ now implies that
$H_\delta$ has ordinals $\bar\kappa<\bar\theta<\kappa$ and a generic elementary
embedding $\sigma\colon H_{\bar\theta}\to H_\theta$ with $\crit\sigma=\bar\kappa$ and
$\sigma(\bar\kappa)=\kappa$. This shows $(iii)$.

\qquad $(iv)\Rightarrow (iii)$: Fix a regular uncountable $\theta>\kappa$ and let
$\delta=(2^{<\theta})^+$. By $(iv)$ there exist ordinals $\bar\kappa<\bar\delta<\kappa$
and a generic elementary embedding $\pi\colon H_{\bar\delta}\to H_\delta$ with
$\crit\pi=\bar\kappa$ and $\pi(\bar\kappa)=\kappa$. Let $\pi(\bar\theta)=\theta$ (we
can assume that $\theta$ is in the range of $\pi$ by taking it to be largest so that
$(2^{<\theta})^+=\delta$). We will argue that $\bar\kappa$ is virtually
$\bar\theta$-supercompact in $H_{\bar\delta}$, so that by elementarity $\kappa$ will be
virtually $\theta$-supercompact in $H_\delta$, and hence also in $V$. Consider the
restriction $\sigma:=\pi\colon H_{\bar\theta}\to H_\theta.$ Note that $H_\theta$ is
closed under ${<}\bar\theta$-sequences (and more) in $V$. We can assume without loss
that $\sigma$ lives in a $\col(\omega,H_{\bar\theta})$-extension. Let $\dot\sigma$ be a
$\col(\omega,H_{\bar\theta})$-name for $\sigma$. Now define

\eq{
X := \bar\theta{+}1 \cup \{x\in H_\theta\mid\exists y\in H_{\bar\theta}\,\exists
p\in\col(\omega, H_{\bar\theta})\, p\forces\dot\sigma(\check y)=\check x\}\in V.
}

Note that $\abs X = \abs{H_{\bar\theta}} = 2^{<\bar\theta}$ and that
$\ran\sigma\subset X$. Now let $\overline X\prec H_\theta$ be such that
$X\subset\overline X$ and $\overline X$ is closed under ${<}\bar\theta$-sequences.
Note that we can find such an $\overline X$ of size $(2^{<\bar\theta})^{<\bar\theta}
= 2^{<\bar\theta}$. Let $\M$ be the transitive collapse of $\overline X$, so that
$\M$ is still closed under ${<}\bar\theta$-sequences and we still have that $\abs\M =
2^{<\bar\theta}<\bar\delta$, making $\M\in H_{\bar\delta}$.

\qquad Countable Embedding Absoluteness~\ref{lemm.ctblabs} then implies that
$H_{\bar\delta}$ has a generic elementary embedding $\sigma^*\colon
H_{\bar\theta}\to\M$ with $\crit\sigma^*=\bar\kappa$ and the proof of Countable
Embedding Absoluteness shows that we can ensure that $\sigma^*(\bar\kappa)>\bar\theta$.
This verifies that $\bar\kappa$ is virtually $\bar\theta$-supercompact in
$H_{\bar\delta}$.
}

\rema[][rema.remarkableequiv]{
  The above proof shows that if $\kappa$ is faintly $(2^{<\theta})^+$-strong, then it
  is virtually $\theta$-supercompact, and if it is virtually
  $(2^{<\theta})^+$-supercompact ala Magidor, then it is virtually
  $\theta$-supercompact. It is open whether they are equivalent level-by-level (see
  Question \ref{ques.remarkableequiv}).
}

There are alternate possible virtualisations of strong cardinals that turn out to be
weaker than our notion. Wilson has proposed a virtualisation of strongness for a
cardinal $\kappa$ defined by the existence of generic embeddings $\sigma:H_\theta\to
\M$ such that $\crit \sigma=\kappa$, $\sigma(\kappa)>\theta$, $H_\theta=
H_\theta^{\M}$, but $\M$ is allowed to be ill-founded. His notion is just like our
notion of faintly strong cardinals but the embeddings can have an ill-founded target.


\qquad Next, we define a virtualisation of the $\alpha$-superstrong cardinals.

\defi{
  Let $\theta$ be a regular uncountable cardinal and $\alpha$ be an ordinal. Then a
  cardinal $\kappa<\theta$ is \textbf{faintly $(\theta,\alpha)$-superstrong} if it is
  faintly $\theta$-measurable, as witnessed by an embedding $\pi:H_\theta\to\N$ with
  $\crit \pi=\kappa$, $H_\theta=H_\theta^{\N}$ and
  $\pi^\alpha(\kappa)\leq\theta$.\footnote{Here we set $\pi^\alpha(\kappa) :=
  \sup_{\xi<\alpha}\pi^\xi(\kappa)$ when $\alpha$ is a limit ordinal.} We replace
  ``faintly'' by \textbf{virtually} when $\N\subset V$, we say that $\kappa$ is
  \textbf{faintly $\alpha$-superstrong} if it is faintly $(\theta,\alpha)$-superstrong
  for \textit{some} $\theta$, and lastly $\kappa$ is simply \textbf{faintly
  superstrong} if it is faintly 1-superstrong.

}

\qquad Recall that a cardinal $\kappa$ is \textbf{virtually rank-into-rank} if there
exists a cardinal $\theta>\kappa$ and a generic elementary embedding $\pi\colon
H_\theta\to H_\theta$ with $\crit\pi=\kappa$.

\prop[N.][prop.omegaSuperstrongRankToRank]{
A cardinal $\kappa$ is virtually $\omega$-superstrong if and only if it is virtually
rank-into-rank.
}
\proof{
Clearly every virtually rank-into-rank cardinal is virtually $\omega$-superstrong by
definition. So suppose that $\kappa$ is virtually $\omega$-superstrong, as witnessed by
a generic elementary embedding $\pi:H_\theta\to \M$ with
$\pi^\omega(\kappa)\leq\theta$, and we let $\lambda = \pi^\omega(\kappa)$. First,
observe that if $a\in H_\lambda$, then $a\in H_{\pi^n(\kappa)}$ for some $n<\omega$,
and hence by elementarity, $\pi(a)\in H_{\pi^{n+1}(\kappa)}\subseteq H_\lambda$. Thus,
the restriction of $\pi$ to $H_\lambda$ maps into $H_\lambda$. We will argue that this
map is elementary. Note that $H_\lambda$ is the union of the elementary chain of the
$H_{\pi^n(\kappa)}$ for $n<\omega$. Thus, $H_\lambda\models\varphi(a)$ implies that
$H_{\pi^n(\kappa)}\models\varphi(a)$ for some $n<\omega$, which implies that
$H_{\pi^{n+1}(\kappa)}\models\varphi(\pi(a))$ by elementarity of $\pi$, which finally
implies that $H_\lambda\models\varphi(\pi(a))$. It follows that the restriction
$\pi^*:H_\lambda\to H_\lambda$ defined by $\pi^*(a)=\pi(a)$ witnesses that $\kappa$ is
virtually rank-into-rank.

}
\prop[N.][prop.superstrong]{
  If $\kappa$ is faintly superstrong, then $H_\kappa$ has a proper class of virtually
  strong cardinals.
}
\proof{
  Fix a regular $\theta>\kappa$ and a generic elementary embedding $\pi\colon
  H_\theta\to\N$ with $\crit\pi=\kappa$, $H_\theta=H_\theta^{\N}$ and
  $\pi(\kappa)\leq\theta$.  Let's argue that $H_{\pi(\kappa)}(=H_{\pi(\kappa)}^{\N})$
  thinks that $\kappa$ is virtually strong.  Fixing $\kappa<\delta<\pi(\kappa)$, we
  have that $\pi\restr H_\delta:H_\delta\to H_{\pi(\delta)}^{\N}$.  In a
  $\col(\omega,H_\delta)$-extension of $\N$, there is an embedding $\pi^*:H_\delta\to
  H_{\pi(\delta)}^{\mathcal N}$ with $\crit \pi^*=\kappa$ and $\pi^*(\kappa)>\delta$.

  \qquad Following the proof of Theorem~\ref{theo.rem}, we can build in $\N$, $X\prec
  H_{\pi(\delta)}^{\mathcal N}$ with $H_\delta\subseteq X$ of size $|H_\delta|$ such
  that letting $\M$ be the collapse of $X$,  we get an embedding $\sigma:H_\delta\to\M$
  with $\crit \sigma=\kappa$, $\sigma(\kappa)>\delta$, $H_\delta\subseteq \M$, and
  $\M\in H_{\pi(\kappa)}$, witnessing that  $\kappa$ is virtually $\delta$-strong in
  $H_{\pi(\kappa)}$. Now since $H_\kappa\prec H_{\pi(\kappa)}$, we have that $H_\kappa$
  thinks that there is a proper class of virtually strong cardinals.

}

The following theorem shows that the existence of ``Kunen inconsistencies'' is
precisely what is stopping pre-strongness from being equivalent to strongness.

\theo[N.][theo.virtchar]{
  Let $\theta$ be a regular uncountable cardinal. Then a cardinal $\kappa<\theta$ is
  virtually $\theta$-pre-strong if and only if one of the following holds.
  \begin{enumerate}
    \item $\kappa$ is virtually $\theta$-strong, or
    \item $\kappa$ is virtually $(\theta,\omega)$-superstrong.
  \end{enumerate}
}
\proof{
  $(\Leftarrow)$ is trivial, so we show $(\Rightarrow)$. Let $\kappa$ be virtually
  $\theta$-pre-strong. Assume $(i)$ fails, meaning that there is a generic elementary
  embedding \hbox{$\pi\colon H_\theta\to\N$} for some transitive $\N\subseteq V$ with
  $H_\theta\subset\N$, $\crit\pi=\kappa$ and $\pi(\kappa)\leq\theta$.

  \qquad First, assume that there is some $n<\omega$ such that $\pi^n(\kappa)=\theta$.
  The proof of Proposition~\ref{prop.superstrong} shows that $\kappa$ is virtually
  strong in $H_{\pi(\kappa)}$. By elementarity, $\kappa$ is virtually strong in
  $H_{\pi(\kappa)}$, and repeating this argument shows that $\kappa$ is virtually
  strong in $H_{\pi^n(\kappa)}=H_\theta$. It also follows, by elementarity, that
  $\pi(\kappa)$ is virtually strong in $H_{\pi^2(\kappa)}$, and by applying
  elementarity repeatedly, we get that $\pi^n(\kappa)=\theta$ is virtually strong in
  $\N$. Note that the condition $\pi^n(\kappa)=\theta$ implies that $\theta$ is
  inaccessible in $\N$. Thus, $H_\theta$ satisfies that there is no largest cardinal,
  and so by elementarity $\N$ does not have a largest cardinal also.

  \qquad Let $\delta=(\theta^+)^{\N}$. In particular, $\theta$ is virtually
  $\delta$-strong in $\N$, and so $\N$ has a generic elementary embedding
  $\sigma:H_\delta^{\N}\to \M$ with $\crit\sigma=\theta$ and $H_\theta\subseteq
  H_\delta^{\N}\subseteq\M$. Thus, $H_\theta\prec H_{\sigma(\theta)}^{\M}$, from which
  it follows that $\kappa$ is virtually strong in $H_{\sigma(\theta)}^{\M}$, and, in
  particular, virtually $\theta$-strong. But $H_{\sigma(\theta)}^{\M}$ must be correct
  about this since $H_\theta^{\M}=H_\theta^{\N}=H_\theta$. But then $\kappa$ is
  actually virtually $\theta$-strong, contradicting our assumption that $(i)$ fails.

  \qquad Next, assume that there is a least $n<\omega$ such that
  $\pi^{n+1}(\kappa)>\theta$. In particular, $\pi^n(\kappa)\leq\theta$. Since
  $\pi(\kappa)\leq\theta$, we have as before that  $\kappa$ is virtually strong in
  $H_{\pi^n(\kappa)}$ and that  $\pi^n(\kappa)$ is virtually strong in
  $H_{\pi^{n+1}(\kappa)}^{\N}$. Since $\pi^{n+1}(\kappa)$ is inaccessible in $\N$,
  $\delta=(\theta^+)^{\N}$ exists. Thus, in $H_{\pi^{n+1}(\kappa)}^{\N}$ there is some
  generic elementary embedding $\sigma:H_\delta^{\N}\to \M$ with $\crit
  \sigma=\pi^n(\kappa)$, $\sigma(\pi^n(\kappa))>\delta$ and $H_\theta\subseteq
  H_\delta^N\subseteq \M$. Thus, by elementarity, we get $H_{\pi^n(\kappa)}  \prec
  H_{\sigma(\pi^n(\kappa))}^{\M}$. Since, as we already argued, $\kappa$ is virtually
  strong in $H_{\pi^n(\kappa)}$ this means that $\kappa$ is also virtually strong in
  $H_{\sigma(\pi^n(\kappa))}^{\M}$, and as $H_\theta^{\M} = H_\theta^{\N} = H_\theta$,
  this means that $\kappa$ is actually virtually $\theta$-strong, contradicting our
  assumption that $(i)$ fails.

 \qquad Finally, assume $\pi^n(\kappa)<\theta$ for all $n<\omega$ and let $\lambda =
 \sup_{n<\omega}\pi^n(\kappa)$. Since $\lambda\leq\theta$, we have that $\kappa$ is
 virtually $(\theta, \omega)$-superstrong by definition.
}

We then get the following consistency result.

\coro[N.][coro.inLfaintlyMeasVirtPreStrong]{
  For any uncountable regular cardinal $\theta$, the existence of a virtually
  $\theta$-strong cardinal is equiconsistent with the existence of a faintly
  $\theta$-measurable cardinal.
}
\proof{
  The above Proposition~\ref{prop.superstrong} and Theorem~\ref{theo.virtchar} show
  that virtually $\theta$-pre-strongs are equiconsistent with virtually
  $\theta$-strongs. Let us now argue that if $\kappa$ is faintly $\theta$-measurable in
  $L$, then $\kappa$ is virtually $\theta$-pre-strong in $L$. Suppose that
  $\pi:L_\theta\to \M$ is a generic elementary embedding with $\M$ transitive and
  $\crit\pi=\kappa$. By elementarity, $\M$ satisfies $V=L$, and hence by absoluteness
  of the construction of $L$ and transitivity of $\M$, $\M=L_\beta$ for some cardinal
  $\beta\geq\theta$. But then trivially we have $L_\theta\subseteq \M$.
}

\coro[N.][coro.kunen-prestrong-strong]{
  The following are equivalent.
  \begin{enumerate}
    \item For every regular uncountable cardinal $\theta$, every virtually
        $\theta$-pre-strong cardinal is virtually $\theta$-strong.
    \item There are no virtually rank-into-rank cardinals.
  \end{enumerate}
}
\proof{
  $(\Leftarrow)$: By Proposition~\ref{prop.omegaSuperstrongRankToRank} being virtually
  $\omega$-superstrong is equivalent to being virtually rank-into-rank.  The above
  Theorem~\ref{theo.virtchar} then implies $(\Leftarrow)$.

  \qquad $(\Rightarrow)$: Here we have to show that if there exists a virtually
  rank-into-rank cardinal, then there exists a $\theta>\kappa$ and a virtually
  $\theta$-pre-strong cardinal which is not virtually $\theta$-strong. Let
  $\bra{\kappa,\theta}$ be the lexicographically least pair such that $\kappa$ is
  virtually rank-into-rank as witnessed by a generic embedding $\pi:H_\theta\to
  H_\theta$, which trivially makes $\kappa$ virtually $\theta$-pre-strong. If $\kappa$
  was also virtually $\theta$-strong, then we would have a generic elementary embedding
  $\pi^*:H_\theta\to \M$ with $\crit \pi^* =\kappa$, $\pi^*(\kappa)>\theta$, and
  $\M\subseteq V$. By Countable Embedding Absoluteness~\ref{lemm.ctblabs}, $\M$ sees
  that $\kappa$ virtually rank-into-rank, but then, using elementarity, this reflects
  down below $\kappa$, showing that the pair $\bra{\kappa,\theta}$ could not have been
  least.
}

\qquad We showed in Theorem~\ref{theo.rem} that faintly strong cardinals and virtually
strong cardinals are equivalent and we showed in
Corollary~\ref{coro.inLfaintlyMeasVirtPreStrong} that, in $L$, faintly measurable
cardinals and virtually pre-strong cardinals are equivalent. As a final result of this
section, we separate the faintly measurable and virtually measurable cardinals. The
separation is trivial in general, as successor cardinals can be faintly measurable and
are never virtually measurable, but the separation still holds true if we rule out this
successor case.

\qquad We also show that a cardinal $\kappa$ may not even be faintly $\kappa^+$-strong,
but at the same time have the property that for every regular $\theta$, there is a
generic embedding $\pi:H_\theta\to \M$ with $\crit \pi=\kappa$, $\pi(\kappa)>\theta$,
and $H_\theta\subseteq \M$. In particular, we don't get that $H_\theta=H_\theta^{\M}$.

\theo[G.][theo.virtualsep]{
  If $\kappa$ is virtually measurable, then there is a forcing extension $V[G]$ in
  which $\kappa$ is inaccessible and faintly measurable, but not virtually measurable.
  If we further assume that $\kappa$ is virtually strong, then, in $V[G]$, for every
  regular $\theta$, there are generic elementary embeddings $\sigma:H_\theta^{V[G]}\to
  \M$ with $\crit \sigma=\kappa$, $\sigma(\kappa)>\theta$, and
  $H_\theta^{V[G]}\subseteq \M$.
}
\proof{
Suppose that $\kappa$ is virtually measurable. This implies, in particular, that for
every regular $\theta>\kappa$, we have generic elementary embeddings
\hbox{$\pi:H_\theta\to \M$} with $\crit \pi=\kappa$ such that $\M\in V$. Thus, by
Proposition~\ref{prop.coll}, we can assume that each generic embedding $\pi$ exists in
a $\col(\omega,H_\theta)$-extension.

\qquad Let $\mathbb P_\kappa$ be the Easton support iteration that adds a Cohen subset
to every regular $\alpha<\kappa$, and let $G\subset\mathbb P_\kappa$ be $V$-generic.
Standard computations show that $\mathbb P_\kappa$ preserves all inaccessible
cardinals. In particular, $\kappa$ remains inaccessible in $V[G]$.

\qquad Fix a regular $\theta\gg\kappa$ and let $h\subseteq \col(\omega,H_\theta)$ be
$V[G]$-generic. In $V[h]$, we must have an elementary embedding $\pi:H_\theta\to \M$
with $\crit \pi=\kappa$ and $\M\in V$, and we can assume without loss that $\M$ is
countable. Obviously, $\pi\in V[G][h]$. Working in $V[G][h]$, we will now lift $\pi$ to
an elementary embedding on $H_\theta[G]=H_\theta^{V[G]}$. To ensure that such a lift
exists, it suffices to find in $V[G][h]$ an $\M$-generic filter for $\pi(\mathbb
P_\kappa)$ containing $\pi''G$.\footnote{This standard lemma is referred to in the
literature as the \textbf{lifting criterion}.} Observe first that $\pi''G=G$ since the
critical point of $\pi$ is $\kappa$ and we can assume that $\mathbb P_\kappa\subseteq
V_\kappa$. Next, observe that $\pi(\mathbb P_\kappa)\cong \mathbb P_\kappa*\mathbb
P_{\text{tail}}$, where $\mathbb P_{\text{tail}}$ is the forcing beyond $\kappa$. Since
$\M[G]$ is countable, we can build an $\M[G]$-generic filter $G_{\text{tail}}$ for
$\mathbb P_\text{tail}$ in $V[G][h]$. Thus, $G*G_{\text{tail}}$ is $\M$-generic for
$\pi(\mathbb P_\kappa)$, and so we can lift $\pi$ to $\pi:H_\theta[G]\to
\M[G][G_{\text{tail}}]$. Since $\theta$ was chosen arbitrarily, we have just shown that
$\kappa$ is faintly measurable in $V[G]$.

\qquad Since generic embeddings witnessing the virtual measurability of $\kappa$ are
$\kappa$-powerset preserving, it suffices to show that we cannot have generic
$\kappa$-powerset preserving embeddings witnessing the faint measurability of $\kappa$
in $V[G]$. Fix a regular $\theta<\bar\theta$ and a generic elementary embedding
$\sigma:H_{\bar\theta}[G]\to \N$ with $\crit\sigma=\kappa$ and
$\p(\kappa)^{V[G]}=\p(\kappa)^{\N}$. By elementarity,
$H_{\sigma(\theta)}^{\N}=\sigma(H_\theta)[\sigma(G)]$ is a forcing extension of
$K=\sigma(H_\theta)$ by $\sigma(G)=G*\bar G_{\text{tail}}\subseteq \sigma(\mathbb
P_\kappa)\cong \mathbb P_\kappa*\bar{\mathbb P}_{\text{tail}}$. Thus, we have the
restrictions $\sigma:H_\theta\to K$ and $\sigma:H_\theta[G]\to K[G][\bar
G_{\text{tail}}]$.

\qquad Let us argue that $\p^{V[G]}(\kappa)\subseteq \p^{K[G]}(\kappa)$, and hence we
have equality. Suppose $A\subseteq\kappa$ in $V[G]$ and let $\dot A$ be a nice $\mathbb
P_\kappa$-name for $A$, which can be coded by a subset of $\kappa$. Since
$\crit\sigma=\kappa$, we have that $\dot A\in K$, and hence $A=\dot A_G\in K[G]$. But
now it follows that the $K[G]$-generic for ${\rm Add}(\kappa,1)$, the forcing at stage
$\kappa$ in $\sigma(\mathbb P_\kappa)$, cannot be in $V[G]$. Thus, we have reached a
contradiction, showing that $\kappa$ cannot be virtually measurable in $V[G]$.

\qquad Now assume further that $\kappa$ is virtually strong. It suffices to simply note
that $G\in\M[G*G_{\text{tail}}]$ so that $H_\theta[G]\subseteq\N[G*G_{\text{tail}}]$ as
well, and since we lifted $\pi$, we still have $\pi(\kappa) > \theta$.
}

\section{Virtual Woodin cardinals and virtual \Vopenka\ Principle}
\label{sect.woodvop}

In this section we will analyse the virtualisations of Woodin cardinals, which can be
seen as ``boldface'' variants of strong cardinals. We also introduce the ``virtual
\Vopenka\ Principle for finite languages'', which it turns out is not equivalent to the
virtual \Vopenka\ principle (see Theorem~\ref{theo.VopenkaForFiniteSeparate}).

\defi{
  Let $\theta$ be a regular uncountable cardinal. Then a cardinal $\kappa<\theta$ is
  \textbf{faintly $(\theta,A)$-strong} for a set $A$ if there is a forcing
  extension containing a transitive set $M$, a set $B$ and an
  elementary embedding $\pi\colon (H_\theta,\in,A\cap H_\theta)\to(\M,\in,B)$ such that
  $\crit\pi=\kappa$, \hbox{$\pi(\kappa)>\theta$}, $H_\theta=H_\theta^{\M}$, and $B\cap
  H_\theta = A\cap H_\theta$. We say that $\kappa$ is \textbf{faintly
  $(\theta,A)$-supercompact} if we further have that ${^{<\theta}\M}\cap V\subset\M$.
}

\defi{
  \label{defi.woodin}
  A cardinal $\delta$ is \textbf{faintly Woodin} if, given any $A\subset H_\delta$,
  there exists a faintly $({<}\delta,A)$-strong cardinal $\kappa$.
}

As before, for both of the above two definitions we substitute ``faintly'' for
\textbf{virtually} when $\M\subset V$, and  ``strong'', ``supercompact'',  and
``Woodin'' for \textbf{pre-strong}, \textbf{pre-supercompact}, and \textbf{pre-Woodin}
when we do not require that $\pi(\kappa)>\theta$.

\defi{
Let $\theta$ be a regular uncountable cardinal. Then a cardinal $\kappa<\theta$ is
\textbf{virtually $(\theta,A)$-extendible} for a set $A$ if there exists a generic
elementary embedding $\pi\colon (H_\theta,\in,A\cap H_\theta)\to(H_\mu,\in,A\cap
H_\mu)$ such that $\crit\pi=\kappa$ and \hbox{$\pi(\kappa)>\theta$}. As usual,  we
substitute ``extendible" for \textbf{pre-extendible} when we do not require that
$\pi(\kappa)>\theta$.
}

\qquad We note in the following proposition that, in analogy with Woodin cardinals,
virtually Woodin cardinals are Mahlo. This property fails for virtually pre-Woodin
cardinals since \cite{Wilson}, together with Theorem~\ref{theo.vopwood} below, shows
that they can be singular.

\prop[Virtualised folklore][prop.woodmahlo]{
  Virtually Woodin cardinals are Mahlo.
}
\proof{
  Let $\delta$ be virtually Woodin. Note that $\delta$ is a limit of weakly compact
  cardinals by Proposition~\ref{prop.virtit}, making $\delta$ a strong limit. As for
  regularity, assume that we have a cofinal increasing function
  $f\colon\alpha\to\delta$ with $f(0)>\alpha$ and $\alpha<\delta$ and note that $f$
  cannot have any closure points. Fix a virtually $({<}\delta,f)$-strong cardinal
  $\kappa<\delta$. We claim that $\kappa$ is a closure point for $f$, which will yield
  the desired contradiction.

\qquad Let $\gamma<\kappa$ and choose a regular $\theta \in (f(\gamma), \delta)$ above
$\kappa$. We then have a generic elementary embedding $\pi\colon (H_\theta,\in,f\cap
H_\theta)\to(\N,\in,f^+)$ with $H_\theta\subset\N$, $\N\subset V$, $\crit\pi=\kappa$,
$\pi(\kappa)>\theta$, and $f^+$ a function such that $f^+\cap H_\theta = f\cap
H_\theta$. But then $f^+(\gamma) = f(\gamma) < \pi(\kappa)$ by our choice of $\theta$,
so elementarity implies that $f(\gamma)<\kappa$, making $\kappa$ a closure point for
$f$, a contradiction. Thus, $\delta$ is inaccessible.

\qquad Next, let us show that $\kappa$ is Mahlo. Let $C\subset\delta$ be a club and let
$\kappa<\delta$ be a virtually $({<}\delta,C)$-strong cardinal. Let $\theta \in (\min
C,\delta)$ be above $\kappa$ and let $\pi\colon (H_\theta,\in,C\cap
H_\theta)\to(\N,\in,C^+)$ be the associated generic elementary embedding having
$C^+\cap \theta=C$. Then for every $\gamma<\kappa$ there exists an element of $C^+$
below $\pi(\kappa)$, namely $\min C$, so by elementarity $\kappa$ is a limit of
elements of $C$, making it an element of $C$. As $\kappa$ is regular, this shows that
$\delta$ is Mahlo.
}

\qquad The well-known equivalence of the ``function definition'' and ``$A$-strong''
definition of Woodin cardinals holds for virtually Woodin cardinals, and the analogue
of the equivalence between virtually strongs and virtually supercompacts allows us to
strengthen this:

\theo[D.-G.-N.][theo.woodin]{
  For an uncountable cardinal $\delta$, the following are equivalent.
  \begin{enumerate}
    \item $\delta$ is virtually Woodin.
    \item For every $A\subset H_\delta$ there exists a virtually
        $({<}\delta,A)$-supercompact $\kappa<\delta$.
    \item For every $A\subset H_\delta$ there exists a virtually
        $({<}\delta,A)$-extendible $\kappa<\delta$.
    \item For every function $f\colon\delta\to\delta$ there are regular cardinals
        $\kappa<\theta<\delta$, such that $\kappa$ is a closure point of $f$, and a
        generic elementary embedding $\pi\colon H_\theta\to\M$ such that
        $\crit\pi=\kappa$, $H_\theta\subset\M$, $\M\subset V$ and
        $\pi(f\restr\kappa)(\kappa)<\theta$.
    \item For every function $f\colon\delta\to\delta$ there are regular cardinals
        $\kappa<\theta<\delta$, such that $\kappa$ is a closure point of $f$, and a
        generic elementary embedding $\pi\colon H_\theta\to\M$ such that
        $\crit\pi=\kappa$, ${^{<\theta}\M}\cap V\subset\M$, $\M\subset V$ and
        $\pi(f\restr\kappa)(\kappa)<\theta$.
    \item For every function $f\colon\delta\to\delta$ there are regular cardinals
        $\bar\theta<\kappa<\theta<\delta$, such that $\kappa$ is a closure point of
        $f$, and a generic elementary embedding $\pi\colon H_{\bar\theta}\to H_\theta$
        with $\pi(\crit\pi)=\kappa$, $f(\crit\pi)<\bar\theta$ and
        $f\restr\kappa\in\ran\pi$.
  \end{enumerate}
}

\begin{figure}
  \begin{center}
    \includegraphics[scale=0.4]{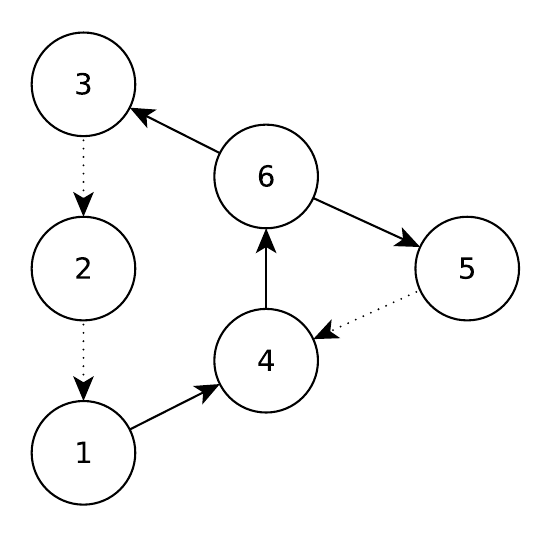}
    \caption{Proof strategy of Theorem~\ref{theo.woodin}, dotted lines are trivial
    implications.}
  \end{center}
\end{figure}

\proof{
  Firstly note that $(iii)\Rightarrow(ii)\Rightarrow(i)$ and $(v)\Rightarrow(iv)$ are
  simply by definition.

  \qquad\framebox{$(i)\Rightarrow(iv)$} Assume $\delta$ is virtually Woodin, and fix a
  function $f\colon\delta\to\delta$. Let $\kappa<\delta$ be virtually $({<}\delta,
  f)$-strong and let $\theta<\delta$ be a regular cardinal such that
  $\sup_{\alpha\leq\kappa}f(\alpha)<\theta$. Then there is a generic elementary
  embedding \hbox{$\pi\colon (H_\theta, \in, f\cap H_\theta) \to (\M, \in, f^+)$} such
  that $H_\theta\subseteq \M$, $f\cap H_\theta=f^+\cap H_\theta$, $\M\subset V$, and
  $\pi(\kappa)>\theta$. Note that, by our choice of $\theta$, $f\restr\kappa\in
  H_\theta$ and $\pi(f\restr\kappa)(\kappa)=f^+(\kappa)=f(\kappa)<\theta$.

 \qquad So it suffices to show that $\kappa$ is a closure point for $f$. Let
 $\alpha<\kappa$. Then
  \eq{
    f(\alpha)=f^+(\alpha)=\pi(f\restr\kappa)(\alpha) =
    \pi(f\restr\kappa)(\pi(\alpha))=\pi(f(\alpha)),
  }

  so $\pi$ fixes $f(\alpha)$ for every $\alpha<\kappa$. Now, if $\kappa$ was not a
  closure point of $f$ then, letting $\alpha<\kappa$ be the least such that
  $f(\alpha)\geq\kappa$, we have
  \eq{
    \theta > f(\alpha) = \pi(f(\alpha)) > \theta,
  }

  a contradiction. Note that we used that $\pi(\kappa)>\theta$ here, so this argument
  would not work if we had only assumed $\delta$ to be virtually pre-Woodin.

  \qquad \framebox{$(iv)\Rightarrow(vi)$} Assume $(iv)$ holds, let
  $f\colon\delta\to\delta$ be given and define $g\colon\delta\to\delta$ as
  $g(\alpha):=(2^{<\gamma_\alpha})^+$, where $\gamma_\alpha$ is the least regular
  cardinal above $|f(\alpha)|$. By $(iv)$ there is a $\kappa<\delta$ which is a closure
  point of $g$ (and so also a closure point of $f$), and there is a regular
  $\lambda\in(\kappa, \delta)$ for which there is a generic elementary embedding
  $\pi\colon H_\lambda\to\M$ with $\crit\pi=\kappa$, $H_\lambda\subset\M$, $\M\subset
  V$, and $\pi(g\restr\kappa)(\kappa)<\lambda$.

  \qquad Let $\theta$ be the least regular cardinal above
  $|\pi(f\restr\kappa)(\kappa)|$, and note that \hbox{$H_\theta\in H_\lambda$} by our
  definition of $g$. Thus, both $H_\theta$ and $H_{\pi(\theta)}^{\M}$ are elements of
  $\M$. An application of Countable Embedding Absoluteness~\ref{lemm.ctblabs} then
  yields that $\M$ has a generic elementary embedding $\pi^*\colon H_\theta^{\M}\to
  H_{\pi(\theta)}^{\M}$ such that $\crit\pi^*=\kappa$, $\pi^*(\kappa)=\pi(\kappa)$,
  $\pi(f\restr\kappa)\in\ran\pi^*$, and $\pi(f\restr\kappa)(\kappa)<\theta$. By
  elementarity of $\pi$, $H_\lambda$ has an ordinal $\bar\theta<\kappa$ and a generic
  elementary embedding $\sigma\colon H_{\bar\theta}\to H_\theta$ with
  $\sigma(\crit\sigma)=\kappa$, $f\restr\kappa\in\ran\sigma$ and
  $f(\crit\sigma)<\bar\theta$, which is what we wanted to show.

  \qquad \framebox{$(vi)\Rightarrow(v)$} Assume $(vi)$ holds and let
  $f\colon\delta\to\delta$ be given. Define $g\colon\delta\to\delta$ as
  $g(\alpha):=\langle(2^{<\gamma_\alpha})^+,f(\alpha)\rangle$, where $\gamma_\alpha$ is
  the least regular cardinal above $|f(\alpha)|$. In particular, $g$ codes $f$. By
  $(vi)$ there exist regular $\bar\kappa<\bar\lambda<\kappa<\lambda$ such that $\kappa$
  is a closure point of $g$ (so also a closure point of $f$) and there exists a generic
  elementary embedding $\pi\colon H_{\bar\lambda}\to H_\lambda$ with $\crit\pi =
  \bar\kappa$, $\pi(\bar\kappa)=\kappa$, $g(\bar\kappa)<\bar\lambda$, and
  $g\restr\kappa\in\ran\pi$. Since $f$ is definable from $g$ and
  $g\restr\kappa\in\ran\pi$, it follows that $f\restr\kappa\in\ran \pi$. So let
  $\pi(\bar f)=f\restr\kappa$ with $\bar f\colon\bar\kappa\to\bar\kappa$. Now observe
  that $\bar f=f\restr\bar\kappa$ since for $\alpha<\bar\kappa$, we have
  $f(\alpha)=\pi(\bar f)(\alpha)=\pi(\bar f(\alpha))=\bar f(\alpha)$.

  \qquad Let $\bar\theta$ be the least regular cardinal above $|f(\bar\kappa)|$. By the
  definition of $g$, we have $H_{\bar\theta}\in H_{\bar\lambda}$. Now, following the
  $(iii)\Rightarrow(ii)$ direction in the proof of Theorem~\ref{theo.rem} we get that
  $H_{\bar\lambda}$ has a generic elementary embedding \hbox{$\sigma\colon
  H_{\bar\theta}\to\M$} with $\M$ closed under ${<}\bar\theta$-sequences from $V$,
  $\crit\sigma=\bar\kappa$, $\sigma(\bar\kappa)> \bar\theta$, and $\sigma(\bar
  f)(\bar\kappa)<\bar\theta$. Let $\pi(\bar\theta)=\theta$ and $\pi(\M)=\N$. Now by
  elementarity of $\pi$, we get that there is a generic elementary embedding
  \hbox{$\sigma^*:H_\theta\to \N$} with $\crit\sigma^*=\kappa$,
  $\sigma^*(\kappa)>\theta$, and $\sigma^*(\pi(\bar
  f))(\kappa)=\sigma^*(f\restr\kappa)(\kappa)<\theta$.

  \qquad \framebox{$(vi)\Rightarrow(iii)$} Let $C$ be the club of all cardinals
  $\alpha$ such that $$(H_\alpha,\in, A\cap H_\alpha)\prec(H_\delta,\in,A).$$ Let
  $f\colon\delta\to\delta$ be given as
  $f(\alpha):=\bra{\gamma_0^\alpha,\gamma_1^\alpha}$, where $\gamma_0^\alpha$ is the
  first limit point of $C$ above $\alpha$ and the $\gamma_1^\alpha$ are chosen such
  that $\{\gamma_1^\alpha\mid\alpha<\beta\}$ encodes $A\cap H_\beta$ for inaccessible
  cardinals $\beta<\delta$.

  \qquad Let $\kappa<\delta$ be a closure point of $f$ such that there are regular
  cardinals \hbox{$\bar\theta<\kappa<\theta$} and a generic elementary embedding
  $\pi\colon H_{\bar\theta}\to H_\theta$ such that $\pi(\crit\pi)=\kappa$,
  $f(\crit\pi)<\bar\theta$, and $f\restr\kappa\in\ran\pi$. Let $\bar\kappa=\crit \pi$.
  The same argument as above gives that $\pi(f\restr \bar\kappa)=f\restr\kappa$. In
  particular, it follows that $\bar\kappa$ is a closure point of $f$, and hence
  $\bar\kappa\in C$. We claim that $\bar\kappa$ is virtually $({<}\delta,
  A)$-extendible. Since $\kappa\in C$ because it is a closure point of $f$, it suffices
  by the definition of $C$ to show that
  \eq{
    (H_\kappa,\in,A\cap H_\kappa)\models\godel{\text{$\bar\kappa$ is virtually $(A\cap
    H_\kappa)$-extendible}}.
  }

  Let $\beta$ be the least element of $C$ above $\bar\kappa$, and note that $\beta$ is
  below $\bar\theta$ since $f(\bar\kappa)<\bar\theta$, and the definition of $f$ says
  that the first coordinate of $f(\bar\kappa)$ is a limit point of $C$ above
  $\bar\kappa$. It then holds that $$(H_{\bar\kappa},\in,A\cap H_{\bar\kappa}) \prec
  (H_\beta,\in,A\cap H_\beta)$$ as both $\bar\kappa$ and $\beta$ are elements of $C$.
  Since $f$ encodes $A$ in the manner previously described and
  $\pi(f\restr\bar\kappa)=f\restr\kappa$, we get that $\pi(A\cap H_{\bar\kappa}) = A\cap
  H_\kappa$, and thus
  \eq{
    (H_\kappa,\in,A\cap H_\kappa) \prec (H_{\pi(\beta)},\in,A^*)
  }

  for $A^* := \pi(A\cap H_\beta)$. Now, as $(H_\gamma,\in, A\cap H_\gamma)$ and
  $(H_{\pi(\gamma)},\in,A^*\cap H_{\pi(\gamma)})$ are elements of $H_{\pi(\beta)}$ for
  every $\gamma<\beta$, Countable Embedding Absoluteness~\ref{lemm.ctblabs} implies
  that $H_{\pi(\beta)}$ sees that $\bar\kappa$ is virtually $({<}\beta,
  A^*)$-extendible. Since

  \eq{
  (H_\beta,\in A\cap H_\beta)\prec (H_{\pi(\beta)},\in A^*),
  }

  it follows that $(H_\beta,\in A\cap H_\beta)$ satisfies that $\bar\kappa$ is
  virtually $({<}\beta, A\cap H_\beta)$-extendible. Finally, since

  \eq{
  (H_\beta,\in, A\cap H_\beta)\prec (H_\kappa,\in,A\cap H_\kappa),
  }

  it follows that $(H_\kappa,\in, A\cap H_\kappa)$ satisfies that $\bar\kappa$ is
  virtually $({<}\kappa,A\cap H_\kappa)$-extendible.
}

As a corollary of the proof, we now have that virtually Woodin cardinals and faintly
Woodin cardinals are equivalent.

\prop[][prop.weakerwoodin]{
A cardinal $\delta$ is virtually Woodin if and only if it is faintly Woodin.
}

Indeed, the argument would work as well if we only assumed existence of generic
embeddings $\pi:(H_\theta,\in,A\cap H_\theta)\to (\M,\in,B)$ such that $\crit
\pi=\kappa$, \hbox{$\pi(\kappa)>\theta$}, $H_\theta\subseteq \M$ and $A\cap
H_\theta=B\cap H_\theta$. In other words, we do not need the a priori stronger
assumption that $H_\theta=H_\theta^{\M}$. Using Proposition~\ref{prop.weakerwoodin}, it
suffices to observe that if $\pi:(H_\theta,\in,A)\to (\M,\in,B)$ is a faintly
$(\theta,A)$-strong embedding such that $A$ codes the sequence of $H_\lambda$ for
$\lambda<\theta$, then $H_\theta=H_\theta^{\M}$.

\qquad We should also observe that if $\delta$ is virtually Woodin, then indeed for
every $A\subseteq H_\delta$, we have stationary many virtually
$({<}\delta,A)$-extendible cardinals by an argument very similar to the proof that
virtually Woodin cardinals are Mahlo.

\qquad Recall that the \textbf{virtual Vop\v enka Principle} states that for every
proper class $\mathcal C$ consisting of structures in a common language, there are
distinct structures $\M,\N\in \mathcal C$ for which there is a generic elementary
embedding $\pi\colon\M\to\N$. The second author and Hamkins showed in
\cite{GitmanHamkins} that the virtual \Vopenka\ principle holds if and only if for
every class $A$ there is a proper class of weakly virtually $A$-extendible cardinals
(in our terminology, these are $({<}\on,A)$-pre-extendible cardinals). It follows from
Theorem~\ref{theo.woodin} that if $\on$ is faintly Woodin, then the virtual \Vopenka\
Principle holds. However, since it is consistent that the virtual \Vopenka\ Principle
holds and $\on$ is not Mahlo \cite{GitmanHamkins}, the two assertions are not
equivalent.

\qquad It turns out, however, that a weakening of the virtual \Vopenka\ Principle is
equivalent to the assertion that $\on$ is faintly pre-Woodin. Our formal setting for
working with classes throughout this article will be the second-order God\"el-Bernays
set theory $\gbc$ whose axioms consist of $\zfc$ together with class axioms consisting
of extensionality for classes, class replacement asserting that every class function
when restricted to a set is a set, global choice asserting that there is a class
well-order of sets, and a weak comprehension scheme asserting that every first-order
formula defines a class.

\defi[][defi.VopenkaFinite]{
The \textbf{virtual Vop\v enka Principle for finite languages} states that for every
proper class $\mathcal C$ consisting of structures in a common \emph{finite} language,
there are distinct structures $\M,\N\in \mathcal C$ for which there is a generic
elementary embedding $\pi\colon\M\to\N$.
}

\qquad The arguments in \cite{GitmanHamkins} show:

\theo[][theo.virtVopFinAext]{
  The virtual \Vopenka\ Principle for finite languages holds if and only if for every
  class $A$, there is a weakly virtually $A$-extendible cardinal.
}
\qquad The \Vopenka\ Principle is equivalent to the \Vopenka\ Principle for finite
languages. Indeed, the \Vopenka\ Principle can be restated in terms of the existence of
elementary embeddings between elements of \textit{natural sequences} \cite{Kanamori}.
But as we will see in the next section, this equivalence relies once again relies on
Kunen's Inconsistency.

\defi[]{
  Say that a class function $f\colon\on\to\on$ is an \textbf{indexing function} if it
  satisfies that $f(\alpha)>\alpha$ and $f(\alpha)\leq f(\beta)$ for all
  $\alpha<\beta$.
}

\defi[]{
  Say that an $\on$-sequence $\vec\M=\bra{\M_\alpha\mid\alpha<\on}$ is \textbf{natural}
  if there exists an indexing function $f^{\vec\M}\colon\on\to\on$ and unary relations
  \hbox{$R^{\vec\M}_\alpha\subset V_{f^{\vec\M}(\alpha)}$} such that $\M_\alpha =
  (V_{f^{\vec\M}(\alpha)}, \in, \{\alpha\}, R_\alpha^{\vec\M})$ for every $\alpha$.
}

\qquad The following theorem shows that the virtual \Vopenka\ Principle for finite
languages holds if and only if $\on$ is virtually pre-Woodin if and only if $\on$ is
faintly pre-Woodin if and only if for every class $A$ there is an weakly virtually
$A$-extendible cardinal. In particular, we get that virtually pre-Woodin and faintly
pre-Woodin cardinals $\delta$ are equivalent, and both are equivalent to the assertion
that for every $A\subseteq H_\delta$, there is a $({<}\delta,A)$-pre-extendible
cardinal. In the next section, in Theorem~\ref{theo.VopenkaForFiniteSeparate}, we will
separate the virtual \Vopenka\ Principle for finite languages from the virtual
\Vopenka\ principle.

\theo[D.-G.-N.][theo.vopwood]{
  The following are equivalent.
  \begin{enumerate}
    \item The virtual Vop\v enka Principle for finite languages holds.
    \item For every class $A$, there is a $({<}\on,A)$-pre-extendible cardinal.
    \item For any natural $\on$-sequence $\vec\M$ there exists a generic elementary
        embedding $\pi\colon\M_\alpha\to\M_\beta$ for some $\alpha<\beta$.
    \item $\on$ is virtually pre-Woodin.
    \item $\on$ is faintly pre-Woodin.
  \end{enumerate}
}
\proof{
By Theorem~\ref{theo.virtVopFinAext}, $(i)$ and $(ii)$ are equivalent.
$(i)\Rightarrow(iii)$ and $(iv)\Rightarrow(v)$ are trivial.

\qquad $(v)\Rightarrow(i)$: Assume that $\on$ is faintly pre-Woodin and fix some class
$\mathcal C$ of structures in a common language. Let $\kappa$ be $({<}\on,\mathcal
C)$-pre-strong. Fix some regular $\theta>\kappa$ such that $H_\theta$ has a structure
$M$ from $\C$ of the $\kappa$-th rank among elements of $\mathcal C$, and fix a generic
elementary embedding $$\pi\colon(H_\theta, \in, \mathcal C\cap H_\theta) \to (\N, \in,
\mathcal C^*)$$ with $\crit \pi=\kappa$, $H_\theta=H_\theta^{\N}$ (note that actually
$H_\theta\subseteq \N$ suffices here), and $\mathcal C\cap H_\theta=\mathcal C^*\cap
H_\theta$.  Now we have that \hbox{$\pi:M\to \pi(M)$} is an elementary embedding (note
that the language is finite and therefore can be assumed to be fixed by $\pi$), and
$\pi(M)\neq M$ because $\pi(M)$ is a structure in $\mathcal C^*$ of rank $\pi(\kappa)$
among elements of $\mathcal C^*$. The structure $(\N,\in,\mathcal C^*)$ believes that
both $M$ and $\pi(M)$ are elements of $\mathcal C^*$, and, by Countable Embedding
Absoluteness~\ref{lemm.ctblabs}, $\N$ has a generic elementary embedding between $M$
and $\pi(M)$. Therefore, $(\N,\in,\mathcal C^*)$ satisfies that there is a generic
elementary embedding between two distinct elements of $\C^*$, and hence
$(H_\theta,\in,\mathcal C\cap H_\theta)$ satisfies that there is a generic elementary
embedding between two distinct elements of $\mathcal C$, and it must be correct about
this.

\qquad $(iii)\Rightarrow(iv)$: Assume $(iii)$ holds and assume that $\on$ is not
virtually pre-Woodin, which means that there exists some class $A$ for which there are
no virtually $({<}\on,A)$-pre-strong cardinals. Define a function $f\colon\on\to\on$ by
setting $f(\alpha)$ to be the least regular $\eta>\alpha$ such that $\alpha$ is not
virtually $(\eta,A)$-pre-strong. Also define $g\colon\on\to\on$ by setting $g(\alpha)$
to be the  least strong limit cardinal above $\alpha$ which is a closure point of $f$.
Note that $g$ is an indexing function, so we can let $\vec\M$ be the natural sequence
induced by $g$ and \hbox{$R_\alpha := A\cap H_{g(\alpha)}$}. $(ii)$ supplies us with
$\alpha<\beta$ and a generic elementary embedding

\eq{
\pi\colon(H_{g(\alpha)},\in,A\cap H_{g(\alpha)})\to (H_{g(\beta)}, \in, A\cap
H_{g(\beta)}).
}

Let us argue that $\pi$ is not the identity map, so that it must have a critical point.
Since $g(\alpha)$ is a closure point of $f$, the structure $(H_{g(\alpha)},\in,A\cap
H_{g(\alpha)})$ can define $f$ correctly. So it must satisfy that $f$ does not have a
strong limit closure point above $\alpha$, but the structure $(H_{g(\beta)}, \in, A\cap
H_{g(\beta)})$ does have such a closure point, namely $g(\alpha)$. Thus, $\pi$ must
have a critical point. Now since $g(\alpha)$ is a closure point of $f$, it holds that
$f(\crit\pi)<g(\alpha)$, so fixing a regular $\theta\in(f(\crit\pi), g(\alpha))$ we get
that $\crit\pi$ is virtually $(\theta, A)$-pre-strong, contradicting the definition of
$f$. Hence $\on$ is virtually pre-Woodin.
}

\section{Virtual Berkeley cardinals}
\label{sect.berkeley}

Berkeley cardinals were introduced by W. Hugh Woodin at the University of California,
Berkeley around 1992, as a large cardinal candidate that would potentially be
inconsistent with \zf. They trivially imply Kunen's Inconsistency and are therefore at
least inconsistent with \zfc, but they have not to date been shown to be inconsistent
with $\zf$. In the virtual setting the virtually Berkeley cardinals, like all the other
virtual large cardinals, are small large cardinals that are downwards absolute to $L$.

\qquad The theorems of this section show that virtually Berkeley cardinals are
precisely the large cardinals which separate virtually pre-Woodin cardinals from
virtually Woodin cardinals analogously to how rank-into-rank cardinals separate
virtually strong cardinals from virtually pre-strong cardinals. To show this, we will
argue that the virtualisation of the notion of the \Vopenka\ filter behaves like the
original one if and only if there are no virtually Berkeley cardinals. We also show
that if the virtual \Vopenka\ Principle holds and $\on$ is not Mahlo, then there is a
virtually Berkeley cardinal. Finally, we use virtually Berkeley cardinals to separate
the virtual \Vopenka\ Principle and the virtual \Vopenka\ principle for finite
languages.

\defi{
   A cardinal $\delta$ is \textbf{virtually proto-Berkeley} if for every transitive set
   $\M$ such that $\delta\subset\M$ there exists a generic elementary embedding
   $\pi\colon\M\to\M$ with $\crit\pi<\delta$.

  \qquad If $\crit\pi$ can be chosen arbitrarily large below $\delta$, then $\delta$ is
  \textbf{virtually Berkeley}, and if $\crit\pi$ can be chosen as an element of any
  club $C\subset\delta$ we say $\delta$ is \textbf{virtually club Berkeley}.
}

Suprisingly, it turns out that the virtually club Berkeley cardinals are precisely the
$\omega$-Erd\"os cardinals. This follows from Lemmata 2.5 and 2.8 in
\cite{RemarkableWilson}.

\theo[Wilson]{
  An $\omega$-Erd\H os is equivalent to a virtually club Berkeley. The least such is
  also the least virtually Berkeley cardinal.\footnote{Note that this also shows that
  virtually club Berkeley cardinals and virtually Berkeley cardinals are
  equiconsistent, which is an open question in the non-virtual context.}
}

\prop{
Virtually (proto)-Berkeley cardinals and virtually club Berkeley cardinals are downward
absolute to $L$. If $0^{\#}$ exists, then every Silver indiscernible is virtually club
Berkeley.
}
\proof{
Downward absoluteness to $L$ follows by  Countable Embedding
Absoluteness~\ref{lemm.ctblabs}. Suppose $0^{\#}$ exists and $\delta$ is a limit Silver
indiscernible. Fix a transitive set $\M\in L$ such that $\delta\subseteq \M$ and a club
$C\subseteq \delta$ in $L$. Let $\lambda\gg\text{rank}(\M)$ be a regular cardinal in
$V$, so that every element of $L_\lambda$ is definable from indiscernibles below
$\lambda$.

\qquad In $V$, we can define a shift of indicernibles embedding $\pi:L_\lambda\to
L_\lambda$ with $\crit\pi\in C$ fixing indiscernibles involved in the definition of
$\M$ (these are easy to avoid since there are finitely many). It follows that
$\pi(\M)=\M$. By Countable Embedding Absoluteness~\ref{lemm.ctblabs}, $L$ must have
such a generic elementary embedding $\sigma:L_\lambda\to L_\lambda$ and it restricts to
$\sigma:\M\to \M$. Thus, $\delta$ is virtually club Berkeley, but then so is every
Silver indiscernible.
}

\qquad Virtually (proto-)Berkeley cardinals turn out to be equivalent to their
``boldface'' versions. The proof is a straightforward virtualisation of Lemma~2.1.12
and Corollary~2.1.13 in \cite{Cutolo}.

\prop[Virtualised Cutolo][prop.boldfaceberkeley]{
  If $\delta$ is virtually proto-Berkeley, then for every transitive set $\M$ such that
  $\delta\subset\M$ and every subset $A\subset\M$ there exists a generic elementary
  embedding $\pi\colon(\M,\in,A)\to(\M,\in,A)$ with $\crit\pi<\delta$. If $\delta$ is
  virtually Berkeley then we can furthermore ensure that $\crit\pi$ is arbitrarily
  large below $\delta$.
}
\proof{
  Let $\M$ be transitive with $\delta\subset\M$ and $A\subset\M$. Let

  \eq{
  \N := \M\cup\{A, \{\{A,x\}\mid x\in \M\}\}\cup \{\{A,x\}\mid x\in \M\}
  }

  and note that $\N$ is transitive. Further, both $A$ and $\M$ are definable in $\N$
  without parameters: the set $A$ is defined as the unique set such that there is a set
  $B$ (namely $\{\{A,x\}\mid x\in \M\}$) all of whose elements are pairs of the form
  $\{A, x\}$, and every set $x$ in $\N$ which does not contain $A$ and which is equal
  to neither $A$ nor $B$, must satisfy that $\{A, x\}\in B$. $M$ is defined in $N$ from
  $A$ as the class containing exactly all elements $x$ such that $A$ is not an element
  of the transitive closure of $x$.

  \qquad But this means that a generic elementary embedding $\pi\colon\N\to\N$ fixes
  both $\M$ and $A$, giving us a generic elementary
  $\sigma\colon(\M,\in,A)\to(\M,\in,A)$ with $\crit\sigma=\crit\pi$, yielding the
  desired conclusion.
}

\coro[][coro.virtBerkeleyImpliesConVirtVopAndOnNotMahlo]{
  If there is a model of $\zfc$ with a virtually Berkeley cardinal, then there is a
  model of $\zfc$ in which the virtual Vop\v enka principle for finite languages holds
  and $\on$ is not Mahlo.
}
\proof{
Using Theorem~\ref{theo.virtuallyBerkeley}, it suffices to show that there is
a model of $\zfc$ with a virtually Berkeley cardinal in which $\on$ is not Mahlo. Take
any model with a virtually Berkeley cardinal, call it $\delta$. If $\on$ is not Mahlo
there, then we are done. Otherwise, $\on$ is Mahlo, so we can let $\kappa$ be the least
inaccessible cardinal above $\delta$. Note that $\delta$ is virtually Berkeley in
$H_\kappa\models\zfc$ and $\on$ is not Mahlo in $H_\kappa$. Thus, $H_\kappa$ is a model
of $\zfc$ with a virtually Berkeley cardinal in which $\on$ is not Mahlo.
}

The following is a straight-forward virtualisation of the usual definition of the Vop\v
enka filter (see e.g. \cite{Kanamori}).

\defi[]{
   The \textbf{virtual Vop\v enka filter} $F$ on $\on$ consists of those classes $X$
   for which there is an associated natural $\on$-sequence $\vec\M^X$ such that
   $\crit\pi\in X$ for any $\alpha<\beta$ and any generic elementary
   $\pi\colon\M_\alpha^X\to\M_\beta^X$.
}

Theorem~\ref{theo.vopwood} shows that the virtual Vop\v enka filter is proper if and
only if the virtual Vop\v enka principle for finite languages holds. The proof of
Proposition~24.14 in \cite{Kanamori} also shows that if the virtual Vop\v enka
principle for finite languages holds, then the virtual Vop\v enka filter is normal.
However, as we will see shortly, unlike the Vop\v enka filter, the virtual Vop\v enka
filter might be proper, but not uniform. This means  that there could be an ordinal
$\delta$ such that every natural sequence of structures has a generic elementary
embedding between two of its structures with critical point below $\delta$.

\qquad Standard proofs of the uniformity of the Vop\v enka filter fail in the virtual
context once again because of the absence of Kunen's Inconsistency. In these arguments,
given an $\on$-length sequence $\bra{\mathcal A_\alpha\mid\alpha<\on}$ of structures in
a common language of size some $\gamma$, we come up with a natural sequences $\bra{
\M_\alpha\mid\alpha<\on}$ such that $\M_\alpha$ codes $\mathcal A_\alpha$. Then given
an elementary embedding \hbox{$\pi:\M_\alpha\to \M_\beta$} with $\crit \pi=\kappa$, we
would like to argue that $\pi$ restricts to an elementary embedding $\pi:\mathcal
A_\alpha\to \mathcal A_\beta$, but this requires that the common language is fixed by
$\pi$. By elementarity of $\pi$, we have that $\mathcal A_\beta$ is a structure in the
language of size $\pi(\gamma)$, and hence $\pi(\gamma)=\gamma$.

\qquad In the presence of Kunen's Inconsistency, it must then be the case that
$\gamma<\kappa$, which implies that the common language is fixed. But no argument like
this will be possible in the virtual context. Indeed, we will show that if the virtual
\Vopenka\ filter is proper, then it is uniform if and only if there are no virtually
Berkeley cardinals if and only if $\on$ is virtually Woodin.

\theo[G.-N.][theo.virtuallyBerkeley]{
  If $\delta$ is a virtually proto-Berkeley cardinal, then for every natural sequence
  $\vec\M=\bra{\M_\alpha\mid\alpha<\on}$, there are ordinals $\alpha<\beta$ and an
  elementary embedding $\pi:\M_\alpha\to \M_\beta$ with $\crit \pi=\kappa<\delta$. In
  particular, the virtual \Vopenka\ principle for finite languages holds and the
  \Vopenka\ filter is proper but not uniform.
}
\proof{
Fix a natural sequence of models $\vec\M=\bra{\M_\alpha\mid\alpha<\on}$. Let
$\theta>\delta$ be a cardinal such that $\vec M\restr\theta\subseteq V_\theta$. Since
$\delta$ is virtually proto-Berkeley, there is an elementary embedding
$\pi:(V_\theta,\in, \vec M\cap V_\theta)\to (V_\theta,\vec M\cap V_\theta)$ with $\crit
\pi=\kappa<\delta$. Then the restriction $\pi:\M_\kappa\to \M_{\pi(\kappa)}$ is an
elementary embedding with critical point $\kappa<\delta$.
}

\coro[][coro.virtBerkeleyImpliesConVirtVopAndOnNotMahlo]{
  If there is a model of $\zfc$ with a virtually Berkeley cardinal, then there is a
  model of $\zfc$ in which the virtual Vop\v enka principle for finite languages holds
  and $\on$ is not Mahlo.
}
\proof{
Using Theorem~\ref{theo.virtuallyBerkeley}, it suffices to show that there is
a model of $\zfc$ with a virtually Berkeley cardinal in which $\on$ is not Mahlo. Take
any model with a virtually Berkeley cardinal, call it $\delta$. If $\on$ is not Mahlo
there, then we are done. Otherwise, $\on$ is Mahlo, so we can let $\kappa$ be the least
inaccessible cardinal above $\delta$. Note that $\delta$ is virtually Berkeley in
$H_\kappa\models\zfc$ and $\on$ is not Mahlo in $H_\kappa$. Thus, $H_\kappa$ is a model
of $\zfc$ with a virtually Berkeley cardinal in which $\on$ is not Mahlo.
}


\lemm[N.][lemm.vopclub]{
 Suppose that the virtual Vop\v enka principle for finite languages holds and that
 there are no virtually Berkeley cardinals. Then the virtually Vop\v enka filter $F$ on
 $\on$ contains every class club $C$.
}
\proof{
  The crucial extra property we get by assuming that there are no virtually Berkeley
  cardinals is that $F$ becomes uniform, i.e., contains every tail
  $(\delta,\on)\subseteq\on$. Indeed, assume that $\delta$ is the least cardinal such
  that $(\delta,\on)\notin F$ (note that $(\gamma,\on)\in F$ up to at least the first
  inaccessible cardinal because a critical point of an elementary embedding
  $\pi:V_\alpha\to V_\beta$ must be inaccessible).

  \qquad Let $\M$ be a transitive set with $\delta\subset \M$ and $\gamma<\delta$ a
  cardinal. As $(\gamma,\on)\in F$ by minimality of $\delta$, we may fix a natural
  sequence $\vec\N$ witnessing this. Let $\vec\M$ be the natural sequence induced by
  the indexing function $f\colon\on\to\on$ given by

  \eq{
    f(\alpha):=\max(\text{rank}(\N_\alpha)+5, \text{rank}(\M)+5)
  }

  and unary relations $R_\alpha := \{\bra{\M, \vec\N_\alpha}\}$, where $\vec
  N_\alpha:=\bra{N_\beta\mid \beta\leq\alpha}$. We can also code in a constant for
  $\delta$ ensuring that any elementary embedding between structures from $\vec M$ must
  fix $\delta$. Suppose $\pi\colon\M_\alpha\to\M_\beta$ is a generic elementary
  embedding with $\crit\pi<\delta$, which exists as $(\delta,\on)\notin F$ and $\pi$
  must fix $\delta$. Since $\pi$ respects the relations $R_\alpha$, we must have
  $\pi\restr\M\colon\M\to\M$ and $\pi(\vec\N_\alpha)=\vec\N_\beta$. We also get that
  $\crit\pi>\gamma$, as

  \eq{
    \pi\restr\N_{\crit\pi}\colon\N_{\crit\pi}\to\N_{\pi(\crit\pi)}
  }

  is an embedding between two structures in $\vec\N$ and $(\gamma,\on)\in F$. This
  means that $\delta$ is virtually Berkeley, a contradiction. Thus $\crit\pi>\delta$,
  implying that $(\delta,\on)\in F$ $\contr$.

  \qquad From here the proof of Lemma~8.11 in \cite{Jech} shows what we wanted.

}

\theo[N.][theo.prewoodwood]{
  If there are no virtually Berkeley cardinals, then $\on$ is virtually pre-Woodin if
  and only if $\on$ is virtually Woodin.
}
\proof{
  Assume $\on$ is virtually pre-Woodin, so the virtual Vop\v enka principle for finite
  languages holds by Theorem~\ref{theo.vopwood}. Let $F$ be the virtual Vop\v enka
  filter, which must be proper. By Lemma~\ref{lemm.vopclub}, $F$ contains every club
  $C$. Assume towards a contradiction that for some class $A$, there are no virtually
  $({<}\on,A)$-extendible cardinals.

  \qquad Define an indexing function $f:\on\to\on$ by $f(\alpha)$ is the least
  $\eta>\alpha$ such that $\alpha$ is not virtually $(\eta,A)$-extendible.  Let $C$ be
  the club of closure points of $f$. Define relations $R_\alpha$ to code $A\cap
  V_{f(\alpha)}$ and $C\cap V_{f(\alpha)}$, and let $\vec M=\bra{
  \M_\alpha\mid\alpha<\on}$ be the associated natural sequence of models. Since $C\in
  F$, there are ordinals $\alpha<\beta$ and an elementary embedding
  \hbox{$\pi:\M_\alpha\to \M_\beta$} with $\crit \pi=\kappa\in C$. It follows that
  $\pi(\kappa)\in C$ as well, and hence it is a closure point of $f$. Thus,
  $f(\kappa)<\pi(\kappa)$. Now consider the restriction

  \eq{
  \pi:(H_{f(\kappa)},A\cap H_{f(\kappa)})\to (H_{f(\pi(\kappa))},\in,A\cap
  H_{f(\pi(\kappa))}),
  }

  which clearly witnesses that $\kappa$ is virtually $(f(\kappa),A)$-extendible,
  contradicting the definition of $f$. Thus, for every class $A$, there is a virtually
  $({<}\on,A)$-extendible cardinal, which implies, by Theorem~\ref{theo.woodin}, that
  $\on$ is virtually Woodin. \footnote{We would like to thank the anonymous referee for
  suggesting this simple proof.}
}

We get the following immediate corollaries from Theorem~\ref{theo.vopwood} and
Proposition~\ref{prop.woodmahlo}:

\coro[][coro.virtVopAndOnNotMahloImpliesVirtBerkeley]{
If the virtual \Vopenka\ Principle for finite languages holds and $\on$ is not Mahlo,
then there is a virtually Berkeley cardinal.
}

\coro{
  The existence of a virtually pre-Woodin cardinal is equiconsistent with the existence
  of a virtually Woodin cardinal.
}

By Corollaries~\ref{coro.virtVopAndOnNotMahloImpliesVirtBerkeley}
and~\ref{coro.virtBerkeleyImpliesConVirtVopAndOnNotMahlo} we then get the following
corollary.

\coro[][coro.virtBerkEquiconsistentVirtVopAndOnNotMahlo]{
  The following are equiconsistent:
  \begin{enumerate}
      \item There is a virtually Berkeley cardinal.
      \item The virtual Vop\v enka principle for finite languages holds and $\on$ is
          not Mahlo.
  \end{enumerate}
}

Next, we observe that even the assumption $\on$ is virtually Woodin is not enough to
guarantee that the virtual \Vopenka\ principle is uniform.

\theo[G.-N.]{
It is consistent that $\on$ is virtually Woodin, but the virtual \Vopenka\ filter is
not uniform.
}
\proof{
Let $V$ be a universe in which $\on$ is virtually Woodin, yet there is a virtually
Berkeley cardinal, for instance, $L$ under the assumption of $0^{\#}$. Then, by
Theorem~\ref{theo.virtuallyBerkeley}, the virtual \Vopenka\ filter cannot be uniform.
}

We would like to thank the anonymous referee for pointing out the following separation
result.

\theo{
It is consistent that the virtual \Vopenka\ Principle holds and $\on$ is Mahlo, but
$\on$ is not virtually Woodin.
}
\proof{
Let $V$ be a universe with a virtually strong cardinal and an $\omega$-Erd\"os cardinal
above it. Assume that $\lambda$ is the least virtually strong cardinal. Since there is
an $\omega$-Erd\"os cardinal above $\lambda$, $V_{\lambda}$ has a proper class of
$\omega$-Erd\"os cardinals, each of which is, in particular, virtually Berkeley.

\qquad To show that the virtual Vop\v enka principle holds in $V_\lambda$, we have to
show that for every $A\subseteq V_\lambda$, there is a proper class of virtually
$({<}\lambda,A)$-pre-strong cardinals. Fix $\gamma<\lambda$ and let
$\gamma<\delta<\lambda$ be a virtually Berkeley cardinal (since these are unbounded in
$\lambda$). We will show that for every $A\subseteq\lambda$, there is a virtually
$({<}\lambda,A)$-pre-strong cardinal $\kappa$ above $\gamma$. Using that $\delta$ is
virtually Berkeley, for every cardinal $\theta\geq\delta$ (below $\lambda$) there
exists a generic elementary embedding

\eq{
\pi_\theta\colon(H_\theta,\in,A\cap H_\theta)\to(H_\theta, \in, A\cap H_\theta)
}

with $\gamma<\crit\pi_\theta<\delta$. By the pigeonhole principle we thus get some
$\gamma<\kappa<\delta$ which is the critical point of unboundedly many $\pi_\theta$
below $\lambda$, showing that $\kappa$ is virtually $({<}\lambda,A)$-pre-strong.

\qquad Thus, the virtual \Vopenka\ Principle holds in $V_\lambda$, and $\lambda$ is at
least weakly compact (since it is a virtual large cardinal), so, in particular, is
Mahlo. But obviously $\on$ cannot be virtually Woodin by the leastness property of
$\lambda$.
}

Finally, we separate the virtual \Vopenka\ Principle for finite languages from the
virtual \Vopenka\ Principle.

\theo[G.-N.][theo.VopenkaForFiniteSeparate]{
It is consistent that the virtual \Vopenka\ Principle for finite languages holds, but
the virtual \Vopenka\ Principle fails.
}
\proof{
Let $V$ be a universe with a virtually Berkeley cardinal $\delta$ and an inaccessible
cardinal above it. Let $\lambda$ be the least inaccessible cardinal above $\delta$. It
is not difficult to see that $\delta$ remains virtually Berkeley in $V_\lambda$, and so
the virtual \Vopenka\ Principle for finite languages holds in $V_\lambda$ by
Theorem~\ref{theo.virtuallyBerkeley}. The virtual \Vopenka\ Principle fails in
$V_\lambda$ because it implies, in particular, that there is a proper class of
inaccessible cardinals.
}

\section{Questions}
\label{sect.questions}

We noted in Remark~\ref{rema.remarkableequiv} that we proved that if $\kappa$ is
faintly $(2^{<\theta})^+$-strong, then it is virtually $\theta$-supercompact, and if it
is virtually $(2^{<\theta})^+$-supercompact ala Magidor, then it is virtually
$\theta$-supercompact. We therefore ask if they are really equivalent for each
$\theta$, or if this kind of ``catching up'' is necessary:

\ques[][ques.remarkableequiv]{
  Are virtually $\theta$-strong cardinals, virtually $\theta$-supercompacts and
  virtually $\theta$-supercompacts ala Magidor all equivalent, for every regular
  $\theta$?
}

We showed in Corollary~\ref{coro.virtBerkEquiconsistentVirtVopAndOnNotMahlo} that a
virtually Berkeley cardinal is equiconsistent with the \Vopenka\ principle for finite
languages and $\on$ not being Mahlo. This naturally leads to the following question,
asking whether this also holds for the virtual \Vopenka\ Principle:

\ques[][ques.vopenkawithmahlo]{
    Does $\con(\zfc + \text{there exists a virtually Berkeley cardinal})$ imply\\
    $\con(\gbc + \text{the virtual \Vopenka\ Principle} + \on\text{ is not Mahlo})$?
}

Question 1.7 in \cite{RemarkableWilson} asks whether the existence of a
non-$\Sigma_2$-reflecting \textit{weakly remarkable} cardinal always implies the
existence of an $\omega$-Erd\H os cardinal. Here a weakly remarkable cardinal is a
rewording of a virtually pre-strong cardinal. Furthermore, Wilson also showed that a
non-$\Sigma_2$-reflecting virtually pre-strong cardinal is equivalent to a virtually
pre-strong cardinal which is not virtually strong. We can therefore reformulate
Wilson's question to the following equivalent question.

\ques[Wilson][ques.wilson]{
  If there exists a virtually pre-strong cardinal which is not virtually strong, is
  there then a virtually Berkeley cardinal?
}

Wilson also showed, in \cite{RemarkableWilson}, that his question has a positive answer
in $L$, which in particular shows that they are equiconsistent. Our results above at
least give a partially positive result:

\coro{
  If for every class $A$ there exists a virtually $A$-pre-strong cardinal, and for some
  class $A$ there is no virtually $A$-strong cardinal, then there exists a virtually
  Berkeley cardinal.
}
\proof{
  The assumption implies by definition that $\on$ is virtually pre-Woodin but not
  virtually Woodin, so Theorem~\ref{theo.prewoodwood} supplies us with the desired
  result.
}

The assumption that there is a virtually $A$-pre-strong cardinal for every class $A$ in
the above corollary may seem a bit strong, but
Theorem~\ref{theo.virtuallyBerkeley} shows that this is necessary, which might
lead one to think that Question~\ref{ques.wilson} could have a negative answer.

\appendix
\section{Chart of virtual large cardinals}

\begin{figure}[H]
  \centering
  \includegraphics[scale = 0.2]{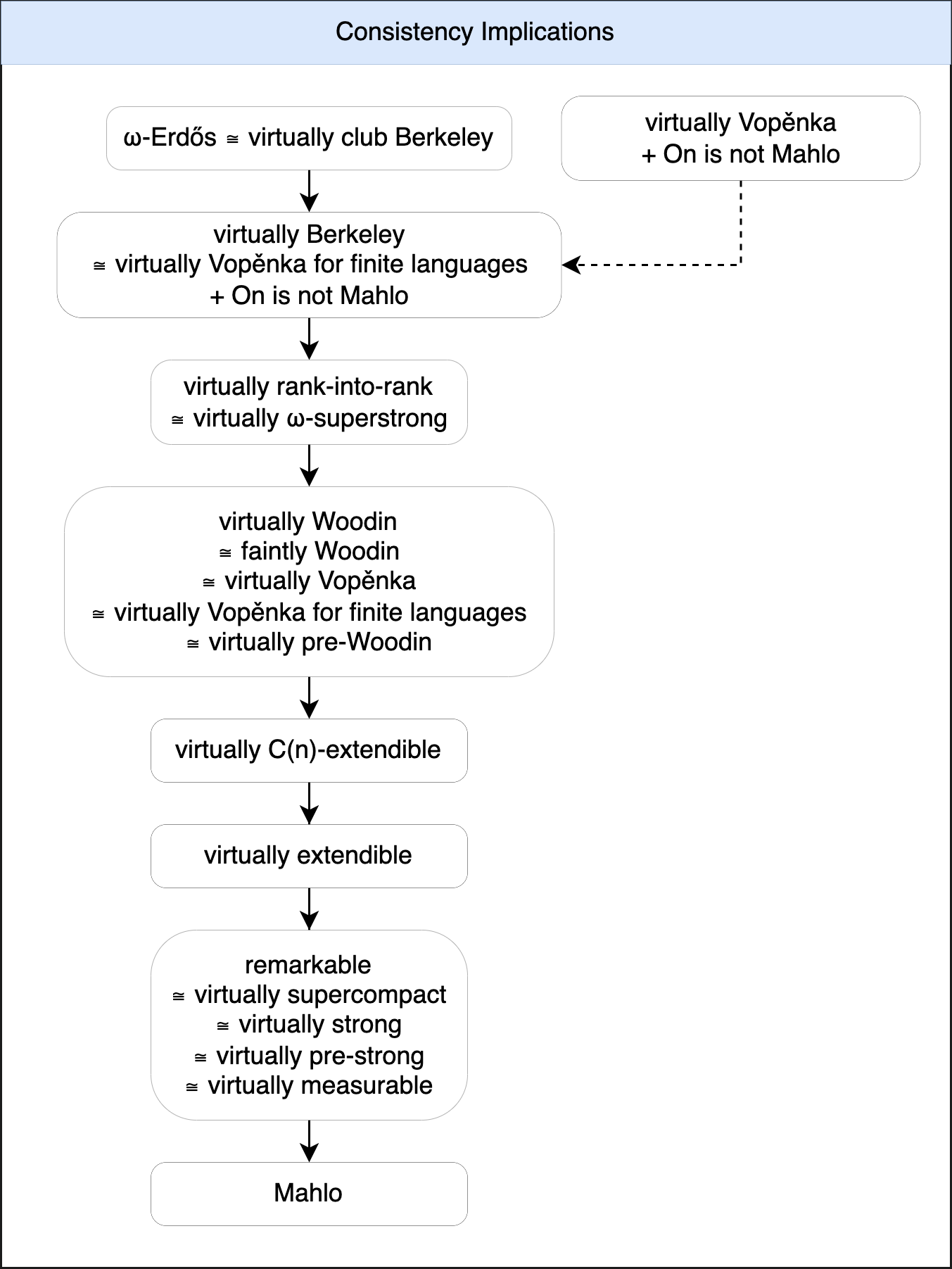}
  \caption{Relative consistency implications between some virtual large cardinals. The
  $\cong$ signs indicate equiconsistency, a solid line indicates that the two are not
equiconsistent, and a dashed line indicates that we do not know whether they are
equiconsistent.}
  \label{fig.consistency-implications}
\end{figure}

\begin{figure}[H]
  \centering
  \includegraphics[scale = 0.2]{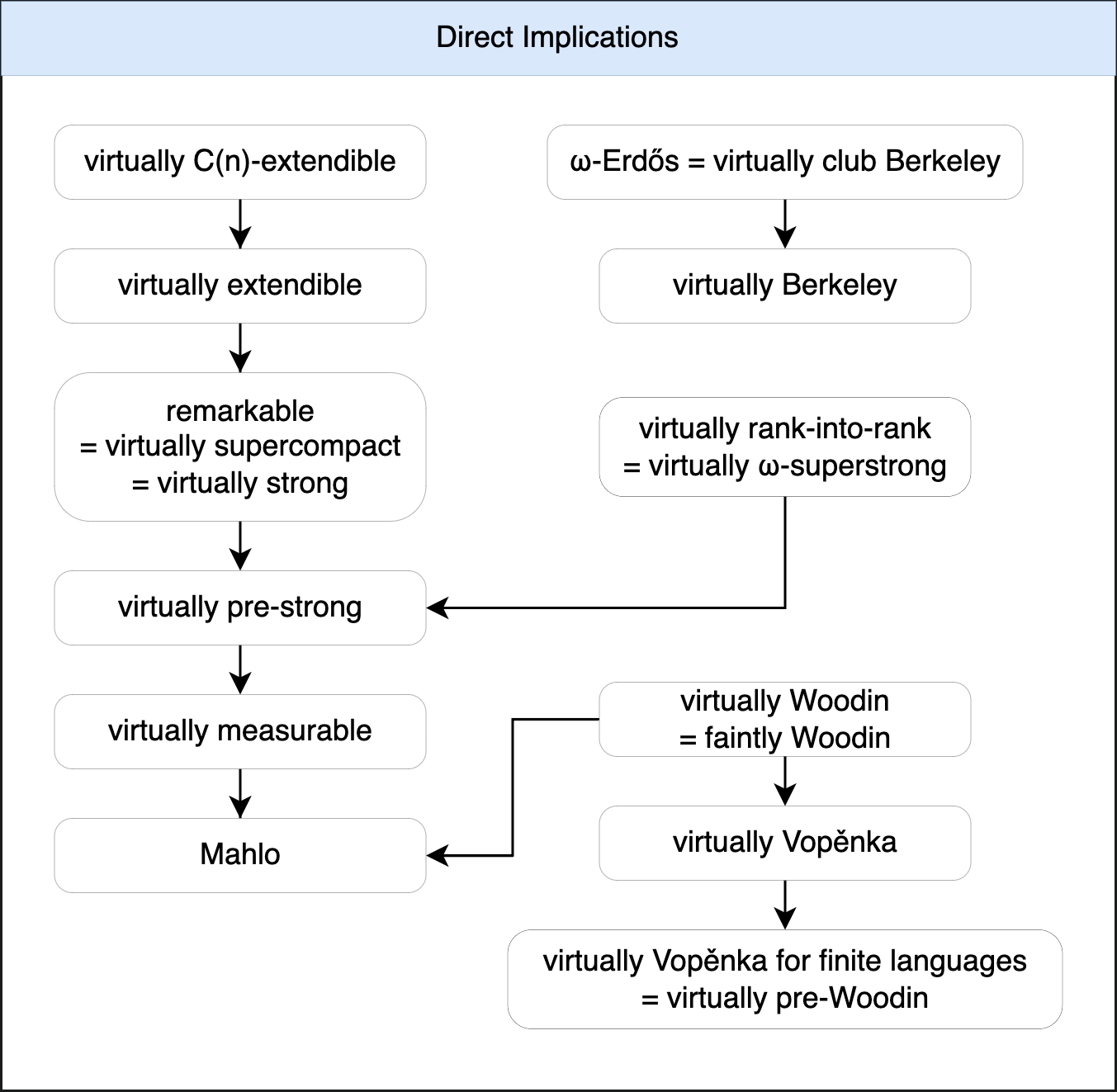}
  \caption{Direct implications between some virtual large cardinals. The equals signs
  indicate equivalence, and a solid line indicates that the two are not equivalent.}
  \label{fig.direct-implications}
\end{figure}

\bibliographystyle{plainnat}
\bibliography{vh}

\end{document}